\documentclass{amsart}

\usepackage{latexsym,  amsfonts, amsmath}

\newtheorem{thm}{Theorem}[section]

\newtheorem{cor}[thm]{Corollary}

\theoremstyle{definition}

\theoremstyle{remark}
\newtheorem{remark}[thm]{Remark}

 \def\Z{{\mathbb{Z}}} \def\N{{\mathbb{N}}}
\def\F{{\mathbb{F}}} \def\G{{\mathbb{G}}} 
\def\C{{\mathbb{C}}}

 \def\GL{{\rm GL}}
\def\E{{\rm E}} \def\D{{\rm D}}
\def\Spec{{\rm Spec}}

\def\ab{{\rm ab}}
\def\deg{{\rm deg}}
\def\End{{\rm End}}
\def\Hom{{\rm Hom}}
\def\M{{\rm M}}

\begin{document}

\sloppy

\title
 {Dieudonn\'{e} Determinants for Skew Polynomial Rings}

\author{Lenny Taelman}
\email{lenny@math.rug.nl}
\address{IWI RuG\\Postbus 800\\9700AV Groningen\\Nederland}




\begin{abstract}
 We observe that the Dieudonn\'{e} determinant induces a non-negative
 degree function on the ring of matrices over a skew polynomial ring.
 We then apply this degree function to two examples. In the first one,
 we find an expression for the rank of the kernel of an endomorphism 
 of the algebraic group $\mathbb{G}_a^n$ over a field of characteristic
 $p>0$. In the second, we calculate the dimension of the solution space
 of linear matrix differential equations.
\end{abstract}

\maketitle

\section{The Dieudonn\'{e} determinant over skew polynomial rings}

In this section we recall some properties of skew polynomial rings and of
matrices over skew fields. We then state and prove the main result.

\subsection{Skew polynomial rings}

This short treatment is based on chapter 1 of \cite{Coh77}. Proofs can
be found in there.

Let $k$ be a field endowed with an endomorphism
$\alpha: k \to k : a \mapsto a^\alpha$ and an
$\alpha$-derivation $\delta$. The latter is an additive map
$\delta:k\to k: a \mapsto a^\delta$
such that 
 $$ (ab)^\delta = a^\delta b^\alpha + ab^\delta. $$
If $\alpha=1$, this is just a derivation.
The \emph{left skew polynomial ring} in $x$ over $k$ associated with
$\alpha$ and $\delta$, denoted by $ R = k[x;\alpha,\delta] $
is the ring whose elements are expressions of the form
 $ a_0 + a_1x + \ldots + a_mx^m $
with $a_i\in k$. Addition is defined in the usual way.
Multiplication is defined by the commutation rule
 $$ xa = a^\alpha x + a^\delta $$
The left euclidean algorithm holds in $R$, relative to the usual degree map
$\deg : R-\{0\} \to \N $. Hence, $R$ is a left principal ideal domain.
When $\alpha$ is an automorphism of the field $k$ then $R$ is also a
right skew polynomial ring and hence also right euclidean.

Since $R$ is a left principal ideal domain, it is left Ore and it has
a unique skew field of fractions $Q$.
It consists of expressions of the form $r/s$ with $s\neq 0$ and
$r/s = r'/s'$ if and only if there exist  non-zero $u,u'$ in
$R$ with $ur=u'r'$ and $us = u's'$. 
The degree map $\deg : Q \to \Z \cup \{\infty\}: r/s \mapsto \deg(r)-\deg(s)$
is well-defined and satisfies
 $$ \deg(tu)=\deg(t)+\deg(u). $$
It factors over $Q \to Q^\ab = Q^\times/[Q^\times,Q^\times] \cup \{0\}$.

\subsection{Matrices over skew polynomial rings}

$\M(n,S)$ is the ring of $n\times n$ matrices over the ring $S$. The subring
of diagonal matrices is denoted by $\D(n,S)= S^n$.
We write $\GL(n,S)$ for the group of invertible elements in $\M(n,S)$ and
$\E(n,S)$ for the subgroup generated by the elementary matrices (corresponding
to the elementary row and column operations).

Consider the ring $\M(n,Q)$ and its subring $\M(n,R)$. On
the former, the Dieudonn\'{e} determinant (\cite{D43}, \cite{Dra83}) defines a
multiplicative function
 $$ \det : \M(n,Q) \to Q^\ab $$
and an exact sequence of groups
 $$ 1 \to \E(n,Q) \to \GL(n,Q) \to (Q^\ab)^\times \to 1 $$
Furthermore, $\E(n,Q)$ equals the commutator subgroup of
$\GL(n,Q)$.\footnotemark
\footnotetext{
 Over \emph{every} skew field, $\E$ is the kernel of $\det$. It equals the
 commutator subgroup of $\GL$ unless $n=2$ and the skew field is the finite
 field of two elements.
}
As the function $\det$ is rational in the entries of the matrix, we cannot
hope for a restriction to $\M(n,R)$ with the same properties as its
commutative, polynomial, counterpart. However, consider the map
 $$ \deg \det: \M(n,Q) \to \Z \cup \{\infty\}. $$
It is surjective and satisfies $\deg\det(AB)=\deg\det(A)+\deg\det(B)$.
We make a little observation on its restriction to $\M(n,R)$.

\begin{thm}\label{thm}
 The restriction of the map $\deg\det$ to $\M(n,R)$
 induces an exact sequence of unitary semigroups
  $$ 1 \to \GL(n,R) \to \M(n,R) \to \N \cup \{\infty\} \to 0. $$
 Moreover, $\deg\det$ is unique in the sense that every semigroup
 homomorphism $f:\M(n,R)\to\Z\cup\{\infty\}$ that coincides with $\deg\det$
 on $\D(n,R)$ equals $\deg\det$.
\end{thm}

\begin{proof}
 We first prove the uniqueness claim and then use it to deduce the first part
 of the theorem.
 
 {\bf Step 1.}
 Assume
 $f:\M(n,R) \to \Z \cup \{\infty\}$ extends to $\M(n,Q)$. As $\deg$ has a
 unique extension from $R$ to $Q$, we find that $f=\deg\det$ on $\D(n,Q)$.
 Since $\M(n,Q)=\D(n,Q)\M(n,R)$, there can at most be one extension of $f$.
 It is a straightforward computation to show that the extension given by
 $f(DA)=f(D)+f(A)$ with $D\in\D(n,Q)$ and $A\in\M(n,R)$ is a well-defined
 semigroup homomorphism. Observe that this unique lift
 $f:\M(n,Q)\to \Z\cup\{\infty\}$ has to factor over $\det:\M(n,Q)\to Q^\ab$.
 Because $\D(n,Q)$ maps surjectively to $Q^\ab$ and $f$ coincides with
 $\deg\det$ on $\D(n,Q)$, we conclude that $f=\deg\det$.

 {\bf Step 2.} 
 Given a matrix $A\in\M(n,R)$, consider the quotient of left $R$-modules
 $R^n/R^nA$. It is a vector space over $k$. Now define
  $$ f:\M(n,R)\to \N\cup\{\infty\}: A \mapsto \dim_k R^n/R^nA. $$
 Because $f(AB)=f(A)+f(B)$ and because $f(D)=\deg\det(D)$ for diagonal
 matrices $D$, we conclude that $f=\deg\det$ on all of $\M(n,R)$.
 But $f$ takes values in $\N\cup\{\infty\}$ and $\ker f=\GL(n,R)$. This
 finishes the proof.
\end{proof}

\section{Applications}

By interpreting the skew polynomial ring, we can apply the main result to two
special cases. The first has a trivial derivation, the second a trivial
field endomorphism.

\subsection{Endomorphisms of $\G_a^n$ in characteristic $p>0$}

Let $k$ be a field of characteristic $p>0$.
Denote by $\G_a$ the additive
(algebraic) group over $k$. Its ring of endomorphisms is known to be the skew
polynomial ring $R=k[\tau;a\mapsto a^p,0]$. Here $a\in k$ acts on $\G_a$
by multiplication with $a$ and $\tau$ acts as the $p$-th power frobenius.
More generally, we have a canonical isomorphism
 $$ \End(\G_a^n) = \M(n,R). $$
Given a $\phi \in \End(\G_a^n)$, its kernel is a closed subgroup scheme
with an $\F_p$-linear structure. We can now determine its size.
\begin{cor}\label{a}
 If $\deg\det \phi < \infty$ then $\ker \phi$ is a closed subgroup
 scheme of $\G_a^n$ which is finite of rank $p^{\deg\det\phi}$ over
  $\Spec(k)$.
\end{cor}
\begin{proof}
 Consider the function $f:\End(\G_a^n)\to p^\N \cup \{\infty\}$
 that associates with $\phi$ the 
 rank of $\ker\phi$ whenever this kernel is finite and $\infty$ when it
 is not. One readily checks that $f$ is multiplicative and that $f$ coincides
 with $p^{\deg\det}$ for diagonal $\phi$. The result now follows from
 theorem \ref{thm}.
\end{proof}

\begin{remark}[Anderson abelian $t$-modules]
In the theory of Anderson $t$-motives (\cite{And86})
one considers 2 types of objects. They roughly look as follows.
The first object is an \emph{abelian $t$-module of dimension $n$} 
and is an action of $\F_p[t]$ on $\G_a^n/k$.
The second object, called a \emph{$t$-motive of dimension $n$ and rank $r$}
is a module over $R \otimes_{\F_p} \F_p[t]$, which has rank $n$ over
$R=k[\tau;a\mapsto a^p]$ and rank $r$ over $k\otimes \F_p[t]$.
The motives are in some sense dual to the modules and one proves that these
objects actually form anti-equivalent categories. After this has
been proven, one defines the \emph{rank} of an abelian $t$-module to be the
rank of its associated $t$-motive.

However, no straightforward description of the rank of an abelian $t$-module,
other than the rank of its associated motive, seemed to be known for $n>1$.
When $n=1$ we are dealing with Drinfel'd modules and a very simple
description in terms of the degree of the endomorphisms
acting on $\G_a$ is known. Corollary \ref{a} shows that $\deg\det$ is the
sought generalization of this degree function and gives the following
desctiption. The rank of the abelian $t$-module given by
 $$ \phi: \F_p[t] \to \End(\G_a^n/k)=\M(n,R), $$
is the constant $r\in\N$ which satisfies
 $$ (\forall f\in \F_p[t])(\deg\det \phi(f)  =  r \deg f). $$
\end{remark}

\subsection{Systems of linear ordinary differential equations}

Details on the algebraic theory of differential equations can be found
in the book \cite{PS03}.

Let $k$ be a differential field of characteristic $0$. This means $k$ is
equipped with an additive map (called derivation) $k\to k: a \mapsto a'$
satisfying $(ab)'=a'b + ab' $. Let $C \subset k$ be the field of constants,
that is the kernel of the derivation. Assume that the derivation is
non-trivial, i.e. $C \neq k$. Examples are $k=\C(x)$ and $k=\C((x))$
with the usual derivation. Now consider the skew polynomial
ring $R=k[\partial;1,']$ with center $C$. It acts $k$-linearly 
on differential field extensions of $k$ by $\partial a = a'$.

A homogeneous matrix differential equation of the form
 $$ A_0 y + A_1 y' + \ldots + A_d y^{(d)} = 0 $$
where $y$ denotes a vector in $k^n$ and the $A_i$ are matrices in $\M(n,k)$
can be rewritten as
 $$ Ay = 0 $$
with $A=\sum A_i \partial^i \in \M(n,R)$. Conversely, every $A$ corresponds to
such a differential equation. Note that we do not demand $A_d$ to be
invertible. 

As in the proof of theorem \ref{thm}, we associate with $A$ the $R$-module
 $$ M:= R^n/R^nA = R^s \oplus M_{\rm tors}  $$
Take $k \subset l$ to be the Picard-Vessiot extension
of $M_{\rm tors}$, or alternatively, take $l$ to be a universal
differential field extension of $k$. When $s=0$, all solutions of the
differential equation exist over $l$. The \emph{contravariant solution space}
$V$ of $M$ is defined to be the $C$-vector space
 $$ V:=\Hom_R(M,l). $$
It is finite-dimensional if and only if $s=0$, and in that case it is dual to
the $C$-space of solutions in $l^n$ of the given differential equation.
On its dimension, we conclude:

\begin{cor}
 With the notations as above, we have
  $$ \dim_C V = \dim_k M = \deg\det A. $$
\end{cor}

\begin{remark}
 Completely analogous results hold for difference and $q$-difference
 equations.
\end{remark}

\end{document}